\documentclass[12pt,a4paper]{article}
\usepackage{amsfonts}
\usepackage{amsmath}
\usepackage{amssymb}
\usepackage[english]{babel}
\usepackage[T1]{fontenc}
\usepackage{graphicx}
\usepackage{subfigure}
\usepackage{epsfig}
\usepackage{latexsym}
\usepackage{amstext}
\usepackage{amsthm}
\usepackage{graphics}
\usepackage{enumerate}
\usepackage{paralist}
\usepackage{fullpage}
\newtheorem{prop}{Proposition}[section]
\newtheorem{rem}[prop]{Remark}
\newtheorem{lem}[prop]{Lemma}
\newtheorem{theo}[prop]{Theorem}
\newtheorem{cor}[prop]{Corollary}

\newcommand{\beq}{\begin{eqnarray}}
\newcommand{\beqq}{\begin{eqnarray*}}
\newcommand{\eeq}{\end{eqnarray}}
\newcommand{\eeqq}{\end{eqnarray*}}

\def\QED{\quad\hbox{\hskip 4pt\vrule width 5pt height 6pt depth 1.5pt}}

\title{A scaling proof for Walsh's Brownian motion extended arc-sine law}
\author{ {\sc Stavros Vakeroudis}\thanks{Department of Mathematics - Probability and Actuarial Sciences group,
  Universit\'{e} Libre de Bruxelles (ULB), CP210, Boulevard du Triomphe, B-1050 Brussels, Belgium. E-mail: stavros.vakeroudis@ulb.ac.be } \
  {\sc and} \
{\sc Marc Yor}\thanks{Laboratoire de Probabilit\'{e}s et Mod\`{e}les
Al\'{e}atoires (LPMA) CNRS : UMR7599,  Universit\'{e} Pierre et Marie
Curie - Paris VI,  Universit\'{e} Paris-Diderot - Paris VII, 4 Place Jussieu, 75252 Paris Cedex 05, France. \
E-mail: yormarc@aol.com}
\thanks{Institut Universitaire de France, Paris, France. } }
\date{\today}
\begin{document}

\maketitle
\begin{abstract}
We present a new proof of the extended arc-sine law related to Walsh's Brownian motion,
known also as Brownian spider. The main argument mimics the scaling property used previously,
in particular by D. Williams \cite{Wil69}, in the 1-dimensional Brownian case,
which can be generalized to the multivariate case.
A discussion concerning the time spent positive by a skew Bessel process is also presented.
\end{abstract}

$\vspace{5pt}$

$\vspace{5pt}$

\section{Introduction}
\renewcommand{\thefootnote}{\fnsymbol{footnote}}

$\left.a\right)$ Recently, some renewed interest has been shown (see e.g. \cite{PPL12}) in the study of the law of the vector
$$\overrightarrow{A_{1}}=\left(\int^{1}_{0} 1_{(W_{s}\in I_{i})}ds ; \ i=1,2,\ldots,n\right) \ ,$$
where $(W_{s})$ denotes a Walsh Brownian motion, also called Brownian spider (see \cite{Wal78} for Walsh's lyrical description)
living on $I=\bigcup^{n}_{i=1} I_{i}$, the union of $n$ half-lines of the plane, meeting at 0.

For the sake of simplicity, we assume $p_{1}=p_{2}=\ldots=p_{n}=1/n$, i.e.: when returning to 0,
Walsh's Brownian motion chooses, loosely speaking, its "new" ray in a uniform way. In fact, excursion theory and/or the computation of
the semi-group of Walsh's Brownian motion (see \cite{BPY89}) allow to define the process rigorously.

Since $\left(d(0,W_{s});s\geq0\right)$, for $d$ the Euclidian distance, is a reflecting Brownian motion,
we denote by $\left(L_{t},t\geq0\right)$ the unique continuous increasing process such that: \\
$\left(d(0,W_{s})-L_{s};s\geq0\right)$ is a $\mathcal{W}_{s}=\sigma\left\{W_{u},u\leq s\right\}$ Brownian motion.
\\
Let $$\overrightarrow{A_{t}}=\left(A^{(1)}_{t},A^{(2)}_{t},\ldots,A^{(n)}_{t}\right)$$
denote the random vector of the times spent in the different rays.
In Section \ref{secmain} we will state and prove our main Theorem concerning the distribution
of $\overrightarrow{A_{t}}$ for a fixed time.
Section \ref{secstable} deals with the general case of stable variables,
First, we recall some known results and then we state and prove our main Theorem.
Finally, Section \ref{con} is devoted to some remarks and comments.
\\ \\
{\noindent $\left.b\right)$ \textbf{Reminder on the arc-sine law:}} \\
A random variable $A$ follows the arc-sine law if it admits the density:
\beq
\frac{1}{\pi \sqrt{x(1-x)}} \ 1_{\left[\right.0,1\left.\right)}(x).
\eeq
Some well known representations of an arc-sine variable are the following:
\beq\label{arcsine}
A\stackrel{(law)}{=}\frac{N^{2}}{N^{2}+\hat{N}^{2}}\stackrel{(law)}{=}\cos^{2}(U)\stackrel{(law)}{=}\frac{T}{T+\hat{T}}
\stackrel{(law)}{=}\frac{1}{1+C^{2}},
\eeq
where $N,\hat{N}\thicksim\mathcal{N}(0,1)$ and are independent, $U$ is uniform on $[0,2\pi]$,
$T$ and $\hat{T}$ stand for two iid stable (1/2) unilateral variables, and $C$ is a standard Cauchy variable. \\
With $(B_{t},t\geq0)$ denoting a real Brownian motion, two well known examples of arc-sine distributed variables are:
$$g_{1}=\sup\{t<1:B_{t}=0\}, \ \ \mathrm{and} \ \ A^{+}_{1}=\int^{1}_{0} ds \ 1_{(B_{s}>0)} \ ,$$
a result that is due to Paul Lévy (see e.g. \cite{Lev39,Lev39a,Yor95}).
\\ \\
$\left.c\right)$ This point gives some motivation for Section \ref{secstable}. From (\ref{arcsine}),
one could think that more general studies of the time spent positive by diffusions
may bring 2 independent gamma variables (this because $N^{2}$ and $\hat{N}^{2}$ are distributed like
two independent gamma variables of parameter 1/2),
or 2 independent stable $(\mu)$ variables. It turns out that it is the
second case which seems to occur more naturally. We devote Section \ref{secstable} to this case.

\section{Main result}\label{secmain}
Our aim is to prove the following:
\begin{theo}\label{maintheo}
The random vectors $\overrightarrow{A_{T}}/T$ for:
\\ \\
\begin{inparaenum}[(i)]
    \item $T=t$; \ \
    \item $T=\alpha^{(j)}_{s}=\inf\{t:A^{(j)}_{t}>s\}$; \ \
    \item $T=\tau_{l}$, the inverse local times,
\end{inparaenum}
\\ \\
have the same distribution. In particular, it is specified by the iid stable (1/2) subordinators:
$$\left(\left(A^{(j)}_{\tau_{l}},l\geq0\right);1\leq j \leq n\right).$$
Hence:
\beq\label{distrA}
\overrightarrow{A_{1}}\stackrel{(law)}{=}\frac{\overrightarrow{A_{\tau_{1}}}}{\tau_{1}} \ ,
\eeq
which yields that:
\beq\label{distrAcons}
\overrightarrow{A_{1}}\stackrel{(law)}{=}\left(\frac{T_{j}}{\sum^{n}_{i=1}T_{i}} \ ; \ j\leq n\right) \ ,
\eeq
where $T_{j}$ are iid, stable (1/2) variables.
\end{theo}
\noindent{The law of the right-hand side of (\ref{distrA}) is easily computed,
and consequently so is its left-hand side. We refer the reader to \cite{BPY89a} for explicit expressions of this law,
which for $n=2$ reduces to the classical arc-sine law.}
\\ \\
\noindent{\textbf{Proof of Theorem \ref{maintheo}.}} \\
$\left.a\right)$ Clearly, $(ii)$ plays a kind of "bridge" between $(i)$ and $(iii)$. \\
$\left.b\right)$ We shall work with $\left(\alpha^{(1)}_{s},s\geq0\right)$, the inverse of $\left(A^{(1)}_{t},t\geq0\right)$.
It is more convenient to use the notation $\left(\alpha^{(+)}_{s},s\geq0\right)$ for $\left(\alpha^{(1)}_{s},s\geq0\right)$.
We then follow the main steps of \cite{Yor95} (Section 3.4, p. 42), which themselves are inspired by Williams \cite{Wil69};
see also Watanabe (Proposition 1 in \cite{Wat95}) and Mc Kean \cite{McK75}. \\
$\left(A^{(j)}_{t}\right)$ denotes the time spent in $I_{j}$, for any $j\neq1$.
Since
$$
\begin{cases}
    A^{(j)}_{\alpha^{(+)}_{1}} = A^{(j)}_{\tau(L_{\alpha^{(+)}_{1}})} \stackrel{(law)}{=} (L_{\alpha^{(+)}_{1}})^{2} A^{(j)}_{\tau_{1}} \ , \\
    \alpha^{(+)}_{1} = 1+\sum_{j} A^{(j)}_{\alpha^{(+)}_{1}} \ , \\
    \ \ \ \mathrm{and} \\
    \mathrm{for} \ \ \mathrm{every} \ \ u,t\geq0, \ \ \left(L_{\alpha^{(+)}_{u}}^{2}<t\right)=\left(u<A^{(1)}_{\tau_{\sqrt{t}}}\right) \ ,
\end{cases}
$$
and invoking the scaling property, we can write jointly for all $j$'s:
\beq\label{joint}
\left(A^{(j)}_{\alpha^{(+)}_{1}},L^{2}_{\alpha^{(+)}_{1}},\alpha^{(+)}_{1}\right) &\stackrel{(law)}{=}&
\left(L^{2}_{\alpha^{(+)}_{1}}A^{(j)}_{\tau_{1}},L^{2}_{\alpha^{(+)}_{1}},1+\sum_{j} L^{2}_{\alpha^{(+)}_{1}}A^{(j)}_{\tau_{1}}\right) \nonumber \\ &\stackrel{(law)}{=}& \left(\frac{A^{(j)}_{\tau_{1}}}{A^{(1)}_{\tau_{1}}},\frac{1}{A^{(1)}_{\tau_{1}}},\frac{\tau_{1}}{A^{(1)}_{\tau_{1}}}\right).
\eeq
Dividing now both sides by $\alpha^{(+)}_{1}$
and remarking that: $\alpha^{(+)}_{1}A^{(1)}_{\tau_{1}}=\tau_{1}$, we deduce:
\beq\label{jointAL2}
\frac{1}{\alpha^{(+)}_{1}}\left(A^{(j)}_{\alpha^{(+)}_{1}},L^{2}_{\alpha^{(+)}_{1}}\right)\stackrel{(law)}{=}
\frac{1}{\tau_{1}}\left(A^{(j)}_{\tau_{1}},1\right).
\eeq
With the help of the scaling Lemma below, we obtain:
\beq\label{E1f}
E\left[1_{(W_{1}\in I_{1})}f(\overrightarrow{A_{1}},L^{2}_{1})\right]&=&
E\left[\frac{1}{\alpha^{(+)}_{1}}f\left(\frac{\overrightarrow{A_{\alpha^{(+)}_{1}}}}{\alpha^{(+)}_{1}},
\frac{L^{2}_{\alpha^{(+)}_{1}}}{\alpha^{(+)}_{1}}\right)\right] \nonumber \\
&\stackrel{\mathrm{from} \; (\ref{joint})}{=}&
E\left[\frac{A^{(1)}_{\tau_{1}}}{\tau_{1}}f\left(\frac{\overrightarrow{A_{\tau_{1}}}}{\tau_{1}},
\frac{1}{\tau_{1}}\right)\right].
\eeq
$I_{1}$ may be replaced by $I_{m}$, for any $m\in\{2,\ldots,n\}$.
Adding the $m$ quantities found in (\ref{E1f}) and remarking that:
\beq\label{sumtau}
\tau_{1}=\sum^{n}_{i=1}A^{(i)}_{\tau_{1}} \ ,
\eeq
we get:
\beqq
E\left[f(\overrightarrow{A_{1}},L^{2}_{1})\right]&=&
E\left[f\left(\frac{\overrightarrow{A_{\tau_{1}}}}{\tau_{1}},\frac{1}{\tau_{1}}\right)\right].
\eeqq
which proves (\ref{distrA}). Note that from (\ref{jointAL2}), the latter also equals:
\beqq
E\left[f\left(\frac{\overrightarrow{A_{\alpha^{(+)}_{1}}}}{\alpha^{(+)}_{1}},\frac{L^{2}_{\alpha^{(+)}_{1}}}{\alpha^{(+)}_{1}}\right)\right].
\eeqq
Equality in law (\ref{distrAcons}) follows now easily. Indeed, we denote by $\boldsymbol\nu$ the Itô measure
of the Brownian spider, and we have:
\beq
\boldsymbol\nu=\frac{1}{n} \sum^{n}_{j=1} \nu_{j} \ ,
\eeq
where $\nu_{j}$ is the canonical image of $\mathbf{n}$, the standard Itô measure of the space of the excursions
of the standard Brownian motion, on the space of the excursions on $I_{j}$. Hence, with $\lambda_{j}, \ j=1,\ldots,n$
denoting positive constants:
\beqq
E\left[\exp\left(-\sum^{n}_{j=1}\lambda_{j}A^{(j)}_{\tau_{1}}\right)\right]&=&
\exp\left(-\frac{1}{n}\sum^{n}_{j=1}\int \nu_{j}(d\varepsilon_{j})(1-e^{-\lambda_{j}\nu_{j}})\right) \\
&=& \exp\left(-\frac{1}{n}\sum^{n}_{j=1}\sqrt{2\lambda_{j}}\right) \ ,
\eeqq
thus:
\beqq
\overrightarrow{A_{\tau_{1}}}=\left(A^{(j)}_{\tau_{1}} \ ; \ j\leq n\right)\stackrel{(law)}{=}\left(\frac{1}{n^{2}}T_{j} \ ; \ j\leq n\right).
\eeqq
The latter, using (\ref{sumtau}) yields:
\beqq
\overrightarrow{A_{1}}=\frac{\overrightarrow{A_{\tau_{1}}}}{\tau_{1}}
=\frac{\overrightarrow{A_{\tau_{1}}}}{\sum_{i=1}^{n} A^{(i)}_{\tau_{1}} }
\stackrel{(law)}{=}\left(\frac{T_{j}}{n^{2}\sum_{i=1}^{n}n^{-2}T_{i}} \ ; j\leq n\right),
\eeqq
which finishes the proof.
\hfill \QED
\\ \\
It now remains to state the scaling Lemma which played a role in (\ref{E1f}), and which we lift from
\cite{Yor95} (Corollary 1, p. 40) in a "reduced" form.
\begin{lem}\emph{\textbf{(Scaling Lemma)}}\label{scalinglem}
Let $U_{t}=\int^{t}_{0} ds \theta_{s}$, with the pair $(W,\theta)$ satisfying:
\beq
\left(W_{ct},\theta_{ct};t\geq0\right)\stackrel{(law)}{=}\left(\sqrt{c}W_{t},\theta_{t};t\geq0\right).
\eeq
Then,
\beq
E\left[F\left(W_{u},u\leq1\right)\theta_{1}\right]=
E\left[\frac{1}{\alpha_{1}}F\left(\frac{1}{\sqrt{\alpha_{1}}}W_{v\alpha_{1}},v\leq1\right)\right],
\eeq
where $\alpha_{t}=\inf\{s: U_{s}>t\}$.
\end{lem}

\section{Stable subordinators}\label{secstable}

\subsection{Reminder and preliminaries on stable variables}
In this Section, we consider $S_{\mu}$ and $S'_{\mu}$ two independent stable variables
with exponent $\mu\in(0,1)$, i.e. for every $\lambda\geq0$, the Laplace transform of $S_{\mu}$ is given by:
\beq
E[\exp(-\lambda S_{\mu})]=\exp(-\lambda^{\mu}).
\eeq
Concerning the law of $S_{\mu}$, there is no simple expression for its density
(except for the case $\mu=1/2$; see e.g. Exercise 4.20 in \cite{ChY12}).
However, we have that, for every $s<1$ (see e.g. \cite{Zol94} or Exercise 4.19 in \cite{ChY12}):
\beq
E[(S_{\mu})^{\mu s}]=\frac{\Gamma(1-s)}{\Gamma(1-\mu s)} \ .
\eeq
We consider now the random variable of the ratio of two $\mu$-stable variables:
\beq
X=\frac{S_{\mu}}{S'_{\mu}} \ .
\eeq
Following e.g. Exercise 4.23 in \cite{ChY12}, we have respectively the following formulas for the Stieltjes and the Mellin transforms of X:
\beq
E\left[\frac{1}{1+sX}\right]=\frac{1}{1+s^{\mu}} \ , \ s\geq0 \ , \label{Stieljes} \\
E\left[X^{s}\right]=\frac{\sin(\pi s)}{\mu \sin(\frac{\pi s}{\mu})} , \ 0<s<\mu \ .
\eeq
Moreover, the density of the random variable $X^{\mu}$ is given by (see e.g. \cite{Zol57,Lam58} or Exercise 4.23 in \cite{ChY12}):
\beq\label{density}
P(X^{\mu} \in dy)= \frac{\sin(\pi \mu)}{\pi \mu} \ \frac{dy}{y^{2}+2y \cos(\pi \mu)+1} \ , \ y\geq0,
\eeq
or equivalently:
\beq\label{ratio}
\left(\frac{S_{\mu}}{S'_{\mu}}\right)^{\mu}= (C_{\mu}|C_{\mu}>0),
\eeq
where, with $C$ denoting a standard Cauchy variable and $U$ a uniform variable in $\left[\right.0,2\pi\left.\right)$,
$$C_{\mu}=\sin(\pi \mu)C-\cos(\pi \mu)\stackrel{(law)}{=}\frac{\sin(\pi \mu-U)}{U} \ .$$

\subsection{The case of 2 stable variables}
\begin{figure}
\centering
    \includegraphics[width=0.60\textwidth]{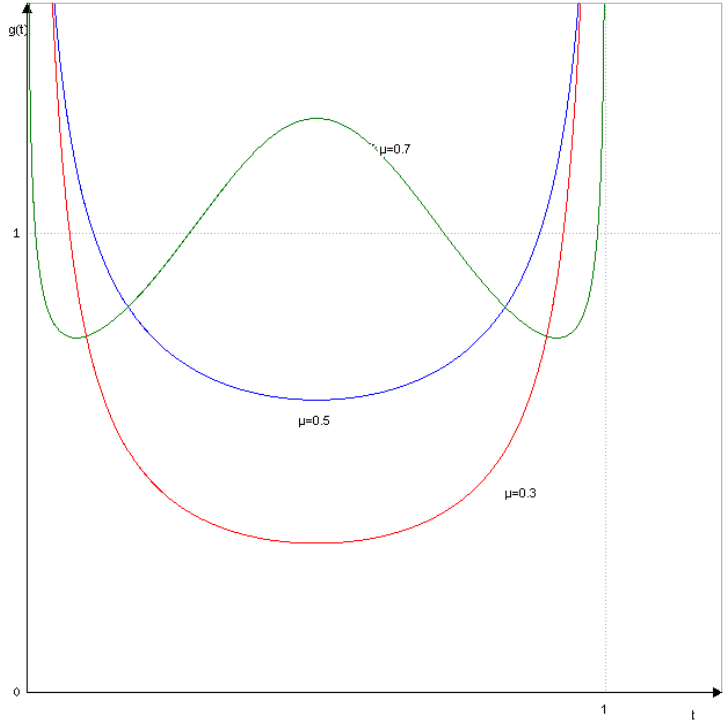}
    \caption{The density function $g$ of $A$, for several values of $\mu$.}
   \label{fig1}
\end{figure}

We turn now our study to the random variable:
\beq\label{A}
A=\frac{S'_{\mu}}{S'_{\mu}+S_{\mu}}=\frac{1}{1+X},
\eeq
\begin{theo}\label{dens}
The density function of the random variable $A$ is given by:
\beq\label{densityfun}
P\left(A\in dz\right)=\frac{\sin(\pi \mu)}{\pi} \frac{dz}{z(1-z)\left[\left(\frac{1-z}{z}\right)^{\mu}+\left(\frac{z}{1-z}\right)^{\mu}+2\cos(\pi\mu)\right]} \ , \ \ z\in[0,1].
\eeq
\end{theo}
\noindent{\textbf{Proof of Theorem \ref{dens}.}} \\
Identity (\ref{A}) is equivalent to:
$$X=\frac{1}{A}-1 \ .$$
Hence, (\ref{Stieljes}) yields:
\beqq
E\left[\frac{1}{1+sX}\right]=E\left[\frac{A}{(1-s)A+s}\right]=\frac{1}{1+s^{\mu}} \ .
\eeqq
We consider now a test function $f$ and invoking the density (\ref{density}) we have $(\nu=\frac{1}{\mu}>1)$:
\beqq
E\left[f\left(\frac{1}{1+X}\right)\right]=
\frac{\sin(\pi \mu)}{\pi \mu} \ \int^{\infty}_{0} \frac{dy}{y^{2}+2y \cos(\pi \mu)+1} \ f\left(\frac{1}{1+y^{\nu}}\right).
\eeqq
Changing the variables $z=\frac{1}{1+y^{\nu}}$, we deduce:
\beqq
E\left[f\left(A\right)\right]=
\frac{\sin(\pi \mu)}{\pi} \ \int^{1}_{0} \frac{dz (1-z)^{\mu-1}}{z^{\mu+1}} \ f\left(z\right) \Delta(z),
\eeqq
where:
\beqq
\Delta(z)&=&\frac{1}{(z^{-1}-1)^{2\mu}+2(z^{-1}-1)^{\mu}\cos(\pi \mu)+1} \\
&=&\frac{z^{2\mu}}{(1-z)^{2\mu}+2(1-z)^{\mu}z^{\mu}\cos(\pi \mu)+z^{2\mu}} \ ,
\eeqq
and (\ref{densityfun}) follows easily.
\hfill \QED
\\ \\
In Figure \ref{fig1}, we have plotted the density function $g$ of $A$, for several values of $\mu$.
\begin{rem}
Similar discussions have been made in \cite{KaY05} in the framework of a skew Bessel process
with dimension $2-2\alpha$ and skewness parameter $p$. Formula (\ref{densityfun}) is a
particular case of formula in \cite{KaY05} for the density of the time spent positive
(called $f_{p,\alpha}$ in \cite{KaY05}).
\end{rem}

\subsection{The case of many stable (1/2) variables}
\begin{figure}
\centering
    \includegraphics[width=0.60\textwidth]{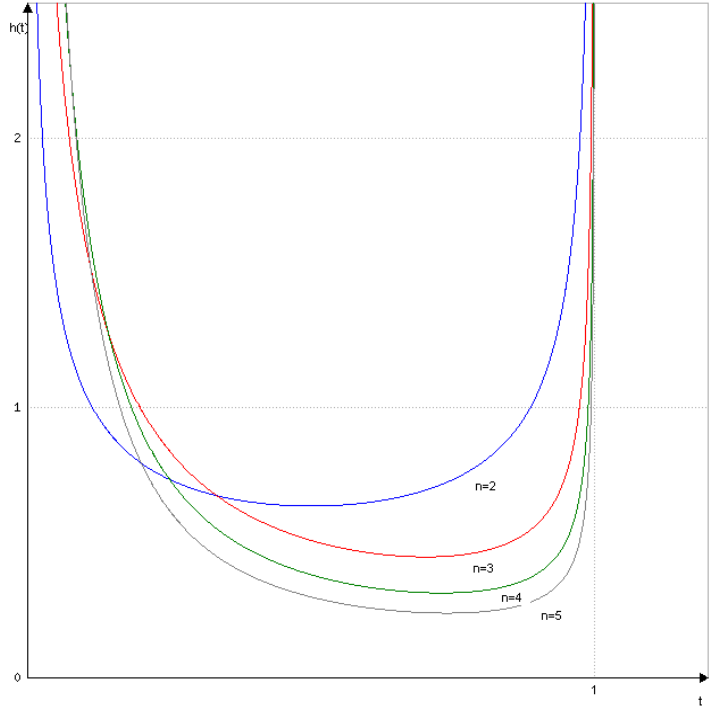}
    \caption{The density function $h$ of $A^{(1)}_{1}$, for several values of $n$.}
   \label{fig2}
\end{figure}

In this Subsection, we consider again $n$ iid stable (1/2) variables, i.e.: $T_{1},\ldots,T_{n}$, and
we will study the distribution of:
\beq
A^{(1)}_{1}=\frac{T_{1}}{T_{1}+\ldots+T_{n}} \ .
\eeq
The following Theorem answers to an open question (and even in a more general sense)
stated at the end of \cite{PPL12}.
\begin{theo}\label{dens2}
The density function of the random variable $A^{(1)}_{1}$ is given by:
\beq\label{densityfun2}
P\left(A^{(1)}_{1}\in dz\right)=\frac{1}{\pi} \frac{dz}{\sqrt{z}\sqrt{1-z}\left[(n-1)z+\frac{1}{n-1}\left(1-z\right)\right]} \ , \ \ z\in[0,1].
\eeq
\end{theo}
\noindent{\textbf{Proof of Theorem \ref{dens2}.}} \\
We first remark that, with $C$ denoting a standard Cauchy variable, using e.g. (\ref{arcsine}):
\beq
A^{(1)}_{1}\stackrel{(law)}{=}\frac{T_{1}}{T_{1}+(n-1)^{2}T_{2}}\stackrel{(law)}{=}\frac{1}{1+(n-1)^{2}C^{2}} \ .
\eeq
Hence, with $f$ standing again for a test function, and invoking the density of a standard Cauchy variable, that is: for every $x\in\mathbb{R}$, $g(x)=\frac{1}{\pi(1+x^{2})}$ we have:
\beqq
E\left[f\left(A^{(1)}_{1}\right)\right]&=&E\left[f\left(\frac{1}{1+(n-1)^{2}C^{2}}\right)\right] \\
&=&\frac{1}{\pi} \ \int^{\infty}_{-\infty} \frac{dx}{1+x^{2}} \ f\left(\frac{1}{1+(n-1)^{2}x^{2}}\right) \\
&\stackrel{x^{2}=y}{=}& \frac{2}{\pi} \ \int^{\infty}_{0} \frac{dy}{2\sqrt{y}(1+y)} \ f\left(\frac{1}{1+(n-1)^{2}y}\right)
\eeqq
Changing the variables $z=\frac{1}{1+(n-1)^{2}y}$, we deduce:
\beqq
E\left[f\left(A^{(1)}_{1}\right)\right]=\frac{1}{\pi} \ \int^{1}_{0} \frac{dz}{(n-1)^{2}z^{2}} \
\frac{(n-1)\sqrt{z}}{\sqrt{z-1}\left(1+\frac{1}{(n-1)^{2}}\left(\frac{1}{z}-1\right)\right)} \ f\left(z\right),
\eeqq
and (\ref{densityfun2}) follows easily.
\hfill \QED
\\ \\
Figure \ref{fig2} presents the plot of the density function $h$ of $A^{(1)}_{1}$, for several values of $n$.

\begin{cor}\label{corA}
The following convergence in law holds:
\beq\label{cvgA}
n^{2}A^{(1)}_{1}(n)\overset{{(law)}}{\underset{n\rightarrow \infty}\longrightarrow} C^{2} \ .
\eeq
\end{cor}
\noindent{\textbf{Proof of Corollary \ref{corA}.}} \\
It follows from Theorem \ref{dens2} by simply remarking that $C\stackrel{(law)}{=}C^{-1}$. Hence:
\beqq
n^{2}A^{(1)}_{1}(n)=\frac{n^{2}}{1+(n-1)^{2}C^{2}}=\frac{1}{\frac{1}{n^{2}}+\left(\frac{n-1}{n}\right)^{2}C^{2}}
\stackrel{n\rightarrow\infty}{\longrightarrow}\frac{1}{C^{2}}\stackrel{(law)}{=}C^{2}.
\eeqq
\hfill \QED

\section{Conclusion and comments}\label{con}
We end up this article with some comments: usually, a scaling argument is "one-dimensional",
as it involves a time-change. Exceptionally (or so it seems to the authors), here we could
apply a scaling argument in a multivariate framework.
We insist that the scaling Lemma plays a key role in our proof.
The curious reader should also look at the totally different proof
of this Theorem in \cite{BPY89a}, which mixes excursion theory and the Feynman-Kac method.

\vspace{30pt}
\noindent\textbf{Acknowledgements} \\
The author S. Vakeroudis is very grateful to Professor R.A. Doney for the invitation
at the University of Manchester as a Post Doc fellow where he prepared a part of this work.



\end{document}